\documentclass{amsart}
\usepackage[utf8]{inputenc} 
\usepackage[T1]{fontenc}    
\usepackage[backref=page]{hyperref}       
\usepackage{amsfonts}       
\usepackage{nicefrac}       
\usepackage{microtype}      

\hypersetup{
  colorlinks   = true, 
  urlcolor     = blue, 
  linkcolor    = teal, 
  citecolor   = teal 
}

\usepackage[usenames,dvipsnames]{xcolor}

\usepackage[color=teal!70!white, bordercolor=teal, textsize=footnotesize,obeyFinal]{todonotes}

\usepackage{marvosym}
\usepackage{amsmath}
\usepackage{amsgen}
\usepackage{amsbsy}
\usepackage{amsopn}
\usepackage{amsthm}
\usepackage{amscd}
\usepackage{amsfonts}
\usepackage{amssymb}
\usepackage{mathrsfs}

\usepackage{ifsym}

\usepackage{tikz-cd}
\usepackage{epigraph}

\usepackage[usestackEOL]{stackengine}
\renewcommand{\epigraphsize}{\small}
\newcommand{\mytextformat}{\epigraphsize\itshape}
\newcommand{\mysourceformat}{\epigraphsize\scshape}

\let\originalepigraph\epigraph 
\renewcommand\epigraph[2]{%
  \setbox0=\hbox{\stackon{\textit{\mytextformat\Longstack{#1}}}%
    {\mysourceformat\scshape\Longstack{#2}}}%
  \ifdim\wd0>.8\linewidth\wd0=.8\linewidth\fi%
  \setlength{\epigraphwidth}{\wd0}%
  \originalepigraph{\textit{#1}}{\textsc{#2}}%
}

\usepackage{color}

\usepackage{graphicx}
\usepackage{comment}

\usepackage{latexsym}

\newtheorem{Theorem}{Theorem}

\newtheorem{Corollary}[Theorem]{Corollary}

\newcommand{\cout}[1]{}

\definecolor{lime}{HTML}{A6CE39}
\DeclareRobustCommand{\orcidicon}{%
    \begin{tikzpicture}
    \draw[lime, fill=lime] (0,0) 
    circle [radius=0.16] 
    node[white] {{\fontfamily{qag}\selectfont \tiny ID}};
    \draw[white, fill=white] (-0.0625,0.095) 
    circle [radius=0.007];
    \end{tikzpicture}
    \hspace{-2mm}
}

\begin{document}


\title[A conservative Turing complete $S^4$ flow]{A conservative Turing complete flow on $S^4$}

\author[P. Su\'arez-Serrato]{P. Su\'arez-Serrato \href{https://orcid.org/0000-0002-1138-0921}{\orcidicon}\\ \today} 

\address{Geometric Intelligence Laboratory, Electrical and Computer Engineering, Univeristy of California, Santa Barbara}

\address{{\it and} Instituto de Matem\'aticas, Universidad Nacional Aut\'onoma de M\'exico UNAM, M\'exico Tenochtitlan}

\begin{abstract}
We present a Turing complete, volume preserving,  smooth flow on the $4$-sphere. The construction relies on the existence of generically rank-2 Poisson structures on $S^4$. A Hamiltonian flow is constructed with the required properties using the symplectic foliation of the Poisson structure. The volume is then preserved because the Poisson structures are unimodular. 
Our main result is a consequence of this construction, showing that such a Turing complete, volume preserving,  smooth flow exists on every compact smooth orientable $4$-manifold. 

\end{abstract} 

\maketitle



The principles behind fluid mechanics date back to the works of ancient philosophers like Aristotle. True advances in our understanding of fluids started being made thanks largely to Leonhard Euler in 1757 \cite{Euler}.
Euler observed that liquid mercury would sink in water whereas solid mercury floated. 
Intuitively, knowing that density must somehow depend on pressure rather than mass alone, Euler set about deriving a formula to calculate density based solely on pressure. 
What resulted was a series of infamously difficult differential equations known today as the Euler Equations.
 With their help, Euler accurately measured the densities of numerous substances, thus laying the foundation for modern day fluid mechanics.
 
Vladimir Igorevich Arnold showed in 1966 that the Euler equations can be regarded as the geodesic flow equations on the group of volume preserving diffeomorphisms, with the right-invariant $L^{2}$ metric defined
by the kinetic energy of the fluid \cite{Arnold66}. 
In this light, the Euler equations describe an evolution in the Lie algebra of divergence-free vector fields, tracing the geodesics in the group volume preserving diffeomorphisms.
This pioneering contribution highlights the importance of the geometry of volume preserving diffeormorphisms and flows for the understanding of Euler's equations, creating the field of topological fluid dynamics.

Nikola Tesla patented a fluid valve in 1920 \cite{Tesla}, and in the 1950s and 1960s the field of fluidics engineered computational machines that brought to life many  of these ideas. 
The applications of microfluidics are now proliferating into multiple scientific and technological domains.
 With these perspectives in mind, one naturally arrives at the question of how general can a fluid computational system be.

The concept of computability has fascinated mathematicians and computer scientists since its inception, leading to many significant breakthroughs in both fields. 
One of the most famous examples of a function that is not computable by standard Turing machines is the Halting Problem \cite{Turing36}. 
It asks whether there exists an algorithm that determines if a given Turing machine will ever halt upon running forever. 
Turing completeness refers to the property of certain computational models, such as Turing machines, that allow them to compute any function that can be expressed by means of a computable procedure. 

A notion of Turing completeness for a flow $\varphi$ on a topological space $X$ was introduced by Moore in 1990 \cite{Moore}.
It requires that for any integer $k>0$, given a Turing machine $T$, an input tape $t$, and a finite string $(t_{1}^{\ast}, \ldots , t_{k}^{\ast})$ of symbols of the alphabet of $T$, there exists an explicitly constructible point $p$ in $X$ and an open set  $U\subset X$ such that the orbit of the flow through $p$ intersects $U$ if, and only if, $T$ halts with an output tape whose positions $-k, \ldots , k$ correspond to the symbols $t_{1}^{\ast}, \ldots , t_{k}^{\ast}$. 

Tao proposed a program to study the relationship between universal computability and the possible existence of singular blowups for Euler's equations of fluid motion  \cite{Tao16, Tao16-2}. 
The overarching idea is to construct a solution that blows up in finite time using controlled, programmable hydrodynamics. 
Tao also gave an example of a coercive potential-well system whose flow is Turing complete \cite{Tao17}. 
 Such dynamics contain orbits for which knowing if they enter a particular set is undecidable. 
 Moreover, explicit descriptions of solutions to the Euler equations with particular dynamics provide evidence for the {\it Turing universality} of such Euler flows \cite{Tao18}.  

Here, we contribute to the understanding of how volume preserving flows can model universal Turing machines, with the next result:

\begin{Theorem}\label{thm:TCvol-pres-flow-S4}
There exists a Turing complete volume preserving smooth flow on $S^4$.
\end{Theorem}

 We induce a volume preserving flow using our previous work that explicitly describes Poisson structures (of generic rank-$2$) on smooth $4$-manifolds \cite{GSV15}. 
 Our construction differs from previously known examples, though it is inspired by the existence of area-preserving Turing complete flows on disks, as demonstrated by Moore.
  It increases the growing list of  volume preserving systems that can perform general computation. 
 Such a result was proven for Euler flows on $T^4$ by Tao \cite{Tao16}, and for every closed oriented $3$-manifold by Cardona, Miranda, Peralta-Salas, and Presas \cite{CMPP21}.

 The main result of this paper, obtained as a consequence of Theorem \ref{thm:TCvol-pres-flow-S4}, is that a compact orientable smooth $4$-manifold also admits such a Turing complete, volume preserving, and Hamiltonian smooth flow:
 
 \begin{Corollary} \label{cor:TCvol-pres-flow-4-mfds}
A closed, oriented, smooth $4$-manifold admits a volume preserving, Turing complete, Hamiltonian smooth flow. 
\end{Corollary}
 
The following paragraphs first review the necessary aspects of Poisson geometry, specifically focusing on unimodularity. 
Next, we will recall how Poisson structures adapted to singular Lefschetz fibrations appear and how we can use them to lift flows to achieve our goals.
After establishing this context, we present our proofs of Theorem \ref{thm:TCvol-pres-flow-S4} and Corollary \ref{cor:TCvol-pres-flow-4-mfds}.

\subsection*{Notions from Poisson geometry}

Write $\frak{X}^{k}$ to denote the smooth $k$-vector fields of the smooth manifold $M$, and $X$ a vector field on $M$. Let $[\![ \cdot \, , \cdot ]\!]$ be the Schouten--Nijenhuis bracket, $[\![ \cdot \, , \cdot ]\!]:  \frak{X}^{k} \times \frak{X}^{l} \to \frak{X}^{k+l - 1}$.
 A bivector field $\pi$ then gives rise to an operator $d_{\pi}: \frak{X}^{\cdot }\to \frak{X}^{\cdot + 1},$
 via the equation $d_{\pi}(X)=[\![ \pi, X]\!]$. Recall that a bivector $\pi$ such that $[\![ \pi, \pi ]\!]=0$ is called a Poisson bivector, and $\pi $ gives rise to a Poisson bracket $\{ \cdot,  \cdot \}_{\pi}$. Informally, a Poisson structure defined by such a bivector furnishes $M$ with a foliation with singularities, and the leaves are symplectic manifolds.
 To a function $h\in C^\infty(M)$ of a Poisson manifold $M$ 
we assign the Hamiltonian vector field $X_h$, defined by $X_h(\cdot )=\{\cdot , h \}_{\pi}$. 
In the simplest case, when the Poisson bracket comes directly from a symplectic form $\omega$, then the entire symplectic manifold is seen as a single leaf, and we have $\{ \cdot , f \} := \omega( \cdot , df)$. 
Thus, we recover the implicit definition of a Hamiltonian flow $X_h$ (with respect to $\omega$), given by the defining equation $\iota_{X_{h}}\omega = -dh$, and both are equivalent.
We have a bundle map $\#_{\pi}:T^{*}M\to TM$ defined by: 
\[
\beta( \#_{\pi} (\alpha)) = \pi(\alpha , \beta), \quad \forall\, \alpha, \beta \in T^{*}M
\]
On the space of differential $1$-forms $\Omega^1(M)$, the tensor $\pi$ induces the following  Lie bracket: 
\[
[ \alpha , \beta ]_{\pi} = L_{\#_{\pi}(\alpha)}\beta - L_{\#_{\pi}(\beta)} \alpha - d( \pi ( \alpha , \beta)),  \quad \forall\, \alpha , \beta \in  \Omega^1(M)
\]

The bundle map $\#_{\pi}$ induces a Lie algebra homomorphism between $[ \cdot , \cdot ]_{\pi}$ and the usual Lie bracket on vector fields.

\subsection*{Unimodular Poisson structures}

This paragraph is inspired by the seminal work of Weinstein on modular fields and classes on Poisson manifolds \cite{Wmod}.
Let ${\mathcal L}$ denote the Lie derivative of the smooth manifold $M$. 
 For a Poisson bivector $\pi$, the equality $\,d_{\pi}^2=0$ holds. Therefore, the {\it Lichnerowicz--Poisson} cochain complex exists, and leads to the definition of the {\it Poisson cohomology spaces} of $(M,\pi)$ \cite{Lich}:
 \[
 H^{k}(M,\pi):=\ker d\pi / {\rm Im}\, d\pi, \quad \, k\in \{ 0, \ldots , \dim M \}
 \]
A vector field $X$ is called a {\it Poisson vector field} if ${\mathcal L}_{X}\pi = 0$. 
Following Weinstein \cite{Wmod}, we can now define the {\it modular vector field}. Let $\Omega$ be a positive volume form for $(M,\pi)$. 
Weinstein showed that, for a smooth function $f$, the operation $f\mapsto {\rm div}_{\Omega}\#_{\pi}(df)$ is a derivation, and hence, a vector field.
 It is called the modular vector field of $(M, \pi )$, with respect to the volume form $\Omega$. 
 This vector field is Poisson, we can denote it by $Y^{\Omega}: C^{\infty}(M)\to C^{\infty}(M)$ and it can also be defined by ${\mathcal L}_{X_f}\Omega = (Y^{\Omega}f)\Omega$. 
 Observe that the first cohomology class of ${\rm div}_{\Omega}\#_{\pi}(df)$ is independent of $\Omega$, because, for a real $a$, 
 \[
 {\rm div}_{a\Omega}\#_{\pi}(d\,\cdot) = {\rm div}_{\Omega}\#_{\pi}(d\,\cdot) + {\rm div}_{a\Omega}\#_{\pi}({\rm log} \,a).
 \]
 
Therefore, the Poisson cohomology class $[Y^{\Omega}]$, called the {\it modular class} of $(M,\pi)$, is a well defined element of the first Poisson cohomology group $H^{1}(M,\pi)$. 
When the modular class vanishes, $[Y^{\Omega}]=0$, the Poisson manifold $(M,\pi)$ is called {\it unimodular} \cite{KS}. 
In an orientable, unimodular, Poisson manifold there exists a volume form that is invariant under all Hamiltonian vector fields \cite{Wmod}.
For example, a symplectic manifold is unimodular, and its Liouville volume form is invariant under every Hamiltonian vector field. 
Przybysz showed that a Poisson structure defined on the algebra of polynomials in $n$ variables which has $(n-2)$ functionally independent Casimirs is unimodular \cite{Przybysz}. 

Consider the $2$-disc $D$, with a symplectic form $\omega_{D}$, and let $(F, \omega_{F})$ be a compact symplectic manifold. 
On the product $D\times F$ let the projections on each factor be
$$p_{D}:D\times F\to D , \,\, p_{F}:D\times F\to F.$$
A symplectic form $\omega_{D\times F}$ on $D\times F$ can be defined by pulling back $\omega_{D}$ and $\omega_{F}$, with respect to $p_{D}$ and $p_{F}$, in that order. 
Define $\omega_{D\times F} := p_{D}^{\ast}\omega_{D}+p_{F}^{\ast}\omega_{F}$.

On the particular case of a product $(D\times F, \omega_{D\times F})$ of a disc $D$ with an orientable surface $F$, consider the bundle given by the projection onto the disc, now renamed $f: D\times F \to D$. On each preimage $f^{-1}(y)=F_{y}$ of a point $y$ in $D$, call $\omega(y)$ the restriction of $\omega_{D\times F}$ to $F_{y}$. Observe that $\omega(y)= \omega_{D\times F}|_{F_{y}} =  \omega_{F}$. Let $\mu$ be an area form of $D$ that agrees with its original orientation. Then, the invariant densities (i.e. the volume forms that are invariant under all Hamiltonian vector fields) of $(D\times F, \omega_{D\times F})$ are of the form 
\begin{equation}\label{eqn:inv-dens}
\omega(y)\wedge f^{*}\mu = \omega_{F} \wedge f^{*}\mu.
\end{equation}

Abouqateb and Boucetta showed that the modular class vanishes on a regular Poisson manifold (one without singularities) when its symplectic foliation admits a compatible Riemannian metric \cite{AB03}.
Further unimodularity criteria around a singular symplectic leaf were shown by Pedroza, Velasco-Barreras, and Vorobiev \cite{PVV18}.

\subsection*{A singular Lefschetz fibration on $S^4$}

We will need a convenient description of $S^4$, which we review now.
Consider the following singular Lefschetz fibration of $f:S^4\to S^2$, discovered by Auroux, Donaldson, and Katzarkov \cite{ADK05}. 
In their construction, the total space of the fibration of $S^4$ is obtained by gluing together three open pieces. 
The first, $X_{-}\simeq T^2\times D^2$, is a trivial fibration that lies over the southern hemisphere base $D_{-}\subset S^2$. 
The second, $X_{+}\simeq S^2\times D^2$, is another trivial fibration that lies over the northern hemisphere base $D_{+}\subset S^2$. 
The third, $W$, lies over a small neighborhood $E$ of the equator on the base $S^2$, such that $W=f^{-1}(E)$. 
In their terminology, the map $f$ has {\it indefinite quadratic singularities} along a 1-dimensional submanifold $\Delta \subset W\subset S^4$. 
This means that around each point of $\Delta$ there exist local co-ordinates $(y_1, y_2, y_3, t)$ on which $\Delta = \{ y_{i}=0 \}$. 
Under $f$, the set $\Delta $ maps to the equator inside $E\subset S^2$.
Moreover, for suitable co-ordinates on the base $S^2$ the map $f$ is represented on $W$ by:
 \begin{equation}\label{eqn:f-on-W}
 (y_1, y_2, y_3, t) \mapsto y_1^2 - (y_2^2 + y_3^2)/2 +it
 \end{equation}
  
 The manifold $W$ is diffeomorphic to a solid torus with a small ball removed. 
 Its boundary components are $\partial_{+}W = T^{2}\times S^1$ and $\partial_{-}W = S^{2}\times S^1$, and $W$ can be seen as a cobordism between these two $3$-manifolds. 
 Over a point $p_{\pm}$ in $f(\partial_{\pm}W )$, denote the fiber $f^{-1}(p_{\pm})$ by $F_{\pm}$. Then $F_{+}\cong T^2$, and $F_{-}\cong S^2$. 
 The transition of the fibers $F_{t}$ lying over a path that joins $p_{+}$ to $p_{-}$ reveals a class of loops that vanish, as $F_{t}$ changes from $F_{+}$ to $F_{-}$. 
 The cycle that manifests this topological change is called the {\it  vanishing cycle}. 
 
 There are multiple possible choices for the diffeomorphisms that ideintify the corresponding boundaries,  $\partial X_{+}$ with $\partial_{+}W$, and $\partial X_{-}$ with $\partial_{-}W$. 
 To obtain $S^4$, two gluing diffeomorphisms need to be specified, one between $X_{-}$ and $W$, and another one between $X_{+}$ and $W$.
  Denote  these gluing diffeomorphisms by $g_{\pm}:X_{\pm}\to W$.
  Attach $X_{+}$ using the non-trivial element in $\pi_1 {\rm Diff}(S^2)$, and to glue $X_{-}$, twist the fibration using a closed path of diffeomorphisms of $T^2$ that corresponds to a unit translation transverse to the vanishing cycle's direction (see \cite{ADK05}). 
Thus put together, $X_{-}\cup W$ is the complement of a closed loop in $S^1\times S^3$. 
The absent loop is isotopic to $S^1\times \{ \ast \} \subset S^1 \times S^3$, projecting non-trivially to the $S^1$-factor.
 Therefore $X_{-}\cup W\simeq S^1\times B^3$, and gluing $X_{+}=D^2\times S^2$ along the boundary as described above yields $S^4$ \cite{ADK05}.

 %
 %

\subsection*{Hamiltonian and the local fibration structure}

Next we state some needed facts about symplectic and Poisson fibrations. 
Consider a locally trivial fiber bundle $p:M\to B$, between closed smooth manifolds $M$ and $B$. Let ${\rm Vert }\subset TM$ be the {\it vertical} sub-bundle formed by the kernels of the differential of $p$:
\[
{\rm Vert }_{x} = \ker d_{x}p = T_{x}F_{p(x)}
\]
Here $F_{b}$ is the fiber of $p$ over $b$. A {\it connection} $\Gamma$ on the fibration $p:M\to B$, in the sense of Ehresmann, means a distribution $\Gamma : x \mapsto {\rm Hor}_{x}$ in $TM$, which splits the tangent bundle,
\[
T_{x}M = {\rm Hor}_{x} \oplus {\rm Vert }_{x}
\]
and obeys the following lifting property. For each point $x$ in $M$ and each base curve $\gamma:[0,1]\to B$ starting at $b_0 = p(x_0)$, there exists an integral curve $\tilde{\gamma}: [0,1]\to  M$ of $\Gamma$ that starts at $x_0$ and covers $\gamma$. A unique {\it horizontal} lift $\gamma_{x}: [0,\varepsilon]\to M$ exists because $T_{x}M$ splits into horizontal and vertical subspaces. %

Observe that every area preserving diffeomorphism $d$ of the $2$-disc $D$ is Hamiltonian.  
Write $H_{\varphi}$ for the Hamiltonian function $H_{\varphi}:D\to {\bf R}$ corresponding to the flow $\varphi$ associated to $d$, so that $\varphi$ is the integral flow of $d H_{\varphi}$.
Let $X$ be the vector field defined by $\iota_{X}\omega_{D}=-d H_{\varphi}$.
Observe that the tangent bundle of $D\times F$ splits as a product, $T(D\times F)\simeq TD\oplus TF$.

Hence we can define a vector field $\overline{X}$ on $D\times S^2$, operates only on the first factor, as follows.
Let $(x,y)\in D\times S^2$, and define  $\bar{\varphi}(x,y):=(\varphi(x), y)$.
Then the Hamiltonian vector field is defined by, 
\begin{equation}\label{eq:Xbar-hamiltonian}
 \\
\iota_{\overline{X}}\omega_{D\times S^2}=-d H_{\bar{\varphi}}.
 \end{equation}
 
 Recall that on a Poisson manifold a Hamiltonian vector field lies tangent to the associated symplectic foliation.  
 
\subsection*{Poisson structures $\pi_{f}$ from singular Lefschetz fibrations $f$}

We will now briefly describe the construction a Poisson structure $\pi_{f}$ on $S^4$ for which the singularities of the fibration $f$ described above coincide with the symplectic foliation singularities of $\pi_{f}$ \cite[Theorem 1.1]{GSV15}. 
These structures are known to be {\it unimodular} \cite{BV21}, a property that will be crucial in our argument. 

Consider the singular circle $\Delta\subset W\subset S^4$, described around equation \eqref{eqn:f-on-W}.  
We have an expression for the Poisson tensor $\pi_{f}$ not only in a neighbourhood of a point on $\Delta$, but on all of normal bundle of $\Delta$.
First we consider the case when the normal bundle of the singular circle is orientable.  
The following expression for $\pi_{f}$  in the normal
bundle of the singular circle inside $W$  was obtained previously \cite[Equation 3.4]{GSV15},
\begin{equation} \label{eq:pi-gamma}
\pi=k(\theta,x_1,x_2,x_3)\left (x_1 \frac{\partial}{\partial x_2} \wedge \frac{\partial}{\partial x_3} + x_2 \frac{\partial}{\partial x_1} \wedge \frac{\partial}{\partial x_3} - x_3 \frac{\partial}{\partial x_1} \wedge \frac{\partial}{\partial x_2} \right ),
\end{equation}
where $k$ is a non-vanishing real function.

Let $N_{\Delta}$ be a neighborhood of $\Delta$, whose normal bundle is orientable. So $N_{\Delta}\cong S^1\times B^3$. On $N_{\Delta}\subset W$, the bivector $\pi_{f}$  is proportional to the product Poisson structure of  $S^1$ with the 
zero Poisson structure, and $B^3$ with the Lie-Poisson structure of $\frak{sl}(2,\mathbb{R})^*$ (relative to an appropriate basis identification \cite{GSV15}).

Second, when the normal bundle of a singular  circle is not orientable, we still have a coordinate description in terms of the quotient by the action of 
\begin{equation}
\label{E:defiota}
i \colon (\theta,x_1,x_2,x_3) \to (\theta+\pi, -x_1, -x_2, x_3).
\end{equation}

The involution $i$ defined on $S^1\times B^3$ by equation \eqref{E:defiota} is a  Poisson map for the bracket  defined by \eqref{eq:pi-gamma}, provided that $k\circ i=k$ \cite[Proposition 3.2]{GSV15}.
The formula for the symplectic form on the symplectic leaf $\Sigma_q$ through 
$q=(\theta,x_1,x_2,x_3)$ at the
point $q$, induced by $\pi$, is (cf \cite[Proposition 3.3]{GSV15}):
\begin{equation}
\label{E:def-omega-Sigma}
\omega_{\Sigma_q}=\frac{1}{k(\theta,x_1,x_2,x_3)\sqrt{x_1^2+x_2^2+x_3^2}}\,\omega_{Area}
\end{equation}
Here, $\omega_{Area}$ is the area form on $\Sigma_q$ induced by the metric $ds^2=d\theta^2+dx_1^2+dx_2^2+dx_3^2$ on $S^1\times B^3$. 

\subsection*{Area preserving diffeomorphisms from generalized shifts}

\begin{figure}
\centering
\includegraphics[width = \textwidth]{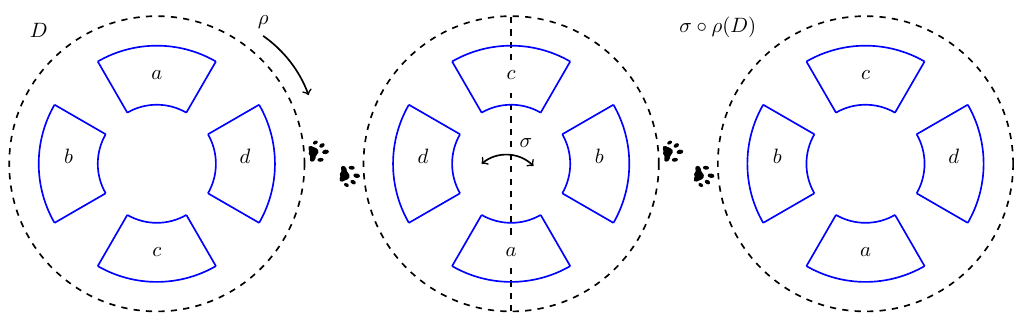}
\caption{Here we illustrate the first couple of steps in a possible parametrization of Moore's diffeomorphism to realize it in an area preserving way. First, with a clockwise $\pi$ rotation $\rho$, followed by exchanging the sectors $b$ and $d$ while preserving the area of small neighborhoods of them with the map $\sigma$. These two transformations, $\rho$ and $\sigma$, all preserve area, so their composition $\sigma\circ \rho$ is also an area preserving diffeomorphism.}
\label{fig:1}
\end{figure}

\begin{figure}
\centering
\includegraphics[width = .7\textwidth]{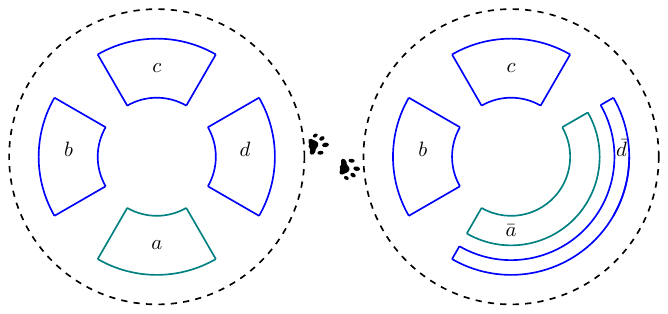}
\caption{Here we show the last step in an area preserving version of Moore's diffeomorphism. The teal mass on the left side sector $a$ circulates to fill the area $\bar{a}$ shown in teal on the right. Likewise, the blue mass on the left side sector $d$ circulates to end occupying the blue area $\bar{d}$ on the lower right}
\label{fig:2}
\end{figure}

Moore showed that any generalized shift is conjugated to a piece-wise linear map, itself made up of finitely many area-preserving components, defined on Cantor blocks, restricted to the Cantor square \cite{Moore}. %
Moore mentioned that such a piece-wise linear map could be extended to a diffeomorphism. 
Cardona, Miranda, Peralta-Salas, and Presas described one such smoothing process, thus unveiling a Turing complete, {\it area-preserving} diffeomorphism $d$ of the $2$-disc $D$ \cite[Theorem 5.2]{CMPP21}. It equals the identity on a neighborhood of the boundary of $D$. Figures (\ref{fig:1}) and (\ref{fig:2}) illustrate the components of Moore's area preserving diffeomorphism $d$, that realizes a generalized shift smoothly.
For brevity, we recommend interested readers consult more details in these two sources we use, by Moore \cite{Moore}, and by Cardona, Miranda, Peralta-Salas, and Presas \cite{CMPP21}.

%
%

\subsection*{Volume preserving and Turing complete flows in dimension 4}

We are now ready to present the proofs of our main results.
 
\begin{proof}[Proof of Theorem \ref{thm:TCvol-pres-flow-S4}]

 Consider the singular Lefschetz fibration of Aroux--Donadson--Katzarkov $f:S^{4}\to S^2$ explained above, assuming the single critical circle $C=f(\Delta)$ lies exactly on the equator of the base $S^{2}$. 
 Let $B$ be a an open set on the base such that $B\subset S^2\setminus \Delta$, such that the fibration trivializes over $B$, so that $f^{-1}(B)\simeq B\times S^2$.
 Here the second factor $S^2$ is the regular fiber of $f$ over $B$.
 Let $D$ be a disk in the fiber $S^{2}\subset f^{-1}(B)\simeq B\times S^2$.
 
  After composing with an area-preserving projection, we may assume that the area-preserving, Turing complete diffeomorphism $d$ described above is acting on $D\subset S^2$. 
  For convenience, we will still call this diffeomorphism $d$.

Denote by $X$ the vector field defined by the derivative of $d$ on $D$. 
Let $\overline{D}$ be a slightly larger closed disc containing $D$. 
Using a bump function interpolating between the boundaries of $\overline{D}$ and $D$ we may control the derivative of $d$ to tends to zero on the boundary as needed so that we can extend $X$ as the zero vector field, from the boundary of $\overline{D}$ to the complement $\overline{D}^{c}$ on the rest of $S^2$. 

The flow $X_{t}$ associated to $X$ preserves the area form $\omega_{S^2}$ of $S^2$.
 Therefore the $1$-form $\iota_{X_{t}}\omega_{S^2}$ is closed for all $t$. 
 Observe that $S^2$ is simply connected, so by Poincar\'e's lemma the form $\iota_{X_{t}}\omega_{S^2}$ is also exact. 
 Hence there exists a function $h_{t}$ such that $dh_{t} = \iota_{X_{t}}\omega_{S^2}$. 
 In other words, both the field $X_{t}$ and the diffeomorphism $d$ are Hamiltonian on $S^2$.

As a consequence of equations (\ref{eqn:inv-dens}) and (\ref{eq:Xbar-hamiltonian}), the vector field $X$ extends and defines a Hamiltonian vector field $\overline{X}$ on $B\times S^2$. 
This vector field $\overline{X}$ may also be extended to all $S^{4}$, setting it equal to the identity on $S^{4}\setminus B\times S^2$, thus defining a flow on all $S^{4}$ which we will still call $\overline{X}$. 
Notice that the flow $\overline{X}$ in $S^4$ only moves along the direction tangent to the symplectic leaves (defined by the fibration structure of $f$). 
Observe that, by construction, the flow  $\overline{X}$ is Hamiltonian with respect to the Poisson structure $\pi_{f}$ defined by $f$ on $S^4$.
Moreover, the Poisson structure $\pi_{f}$ is unimodular.
Therefore there exists a volume form that is invariant under $\overline{X}$ on $S^{4}$.

We will now show the existence of a constructible point and a constructible set to verify the Turing completeness property of the flow $\overline{X}$.
Let $p$ be the constructible point, together with $\mathcal{O}$ its associated constructible open set, guaranteed to exist by the Turing completeness of $X$ on $D\subset S^2$.
Consider a parametrization of $B$ by polar coordinates, and define $c$ to be the center of $B$ in these coordinates. 
Then the point $\bar{p}:= (x,c)$ is a constructible point, for which the open set $\overline{\mathcal{O}}:= B\times  \mathcal{O}$ is a constructible set which together imply that the flow $\overline{X}$ is Turing complete.
The undecidability of the orbit of $p$ under $X$ entering $\mathcal{O}$ implies, by construction, that it is undecidable if the orbit of $\bar{p}$ under $\overline{X}$ enters $\overline{\mathcal{O}}$.

Therefore, we have identified both a constructible point $\bar{p}$ and its corresponding open set $\overline{\mathcal{O}}$, that together verify the Turing completeness of $\overline{X}$.
This concludes the proof of Theorem \ref{thm:TCvol-pres-flow-S4}, as $\overline{X}$ is a Turing complete, Hamiltonian volume preserving flow of $S^4$. \end{proof}

We now present a proof of Corollary \ref{cor:TCvol-pres-flow-4-mfds}:

\begin{proof}[Proof of Corollary \ref{cor:TCvol-pres-flow-4-mfds}]
We mention briefly two possible proofs. 
Let $M$ be a closed oriented smooth $4$-manifold.
First, conseider adding to $M$ a connected sum with $S^4$, where the attaching disc is in the complement of the region of $S^4$ involved in the definition of the flow constructed in Theorem \ref{thm:TCvol-pres-flow-S4}.
 Then, the flow acts on the attached $S^4$, in a way that it can be considered a flow on $M$, extending it as the identity everywhere outside the region originally used in \ref{thm:TCvol-pres-flow-S4}. 

Another proof repeats the construction as in Theorem \ref{thm:TCvol-pres-flow-S4}, except using a generically rank-$2$ Poisson structure $\pi_{f}$ adapted to a singular Lefschetz fibration $f$ on $M$ \cite{GSV15, ST16}. 

In this general case, notice that it is always possible to change a singular Lefschetz fibration into a kind of fibration that is called wrinkled, and which has additional types of singularities. 
These wrinkled fibrations admit Poisson structures that are compatible in the sense that the symplectic leaves are the fibers, and both structures have exactly the same singularities \cite{ST16}.
As before, explicit normal forms for $\pi_{f}$ on neighbourhoods of Lefschetz singularities are available \cite{ST16}.

After performing wrinkled fibration moves, we may assume we have a new fibration $\bar{f}$ from $X$ to $S^2$, which has an open region where the fibre is $S^2$.
Choose a disc region $B$ on the base of the fibration $\bar{f}$, so that the preimage $f^{-1}(B)$ is a product $B\times S^2$.  
Then the arguments in Theorem \ref{thm:TCvol-pres-flow-S4} carry through exactly the same.
Specifically, by Equation \eqref{eq:Xbar-hamiltonian} and the extension arguments mentioned previously, there exists a Turing complete Hamiltonian flow with respect to the Poisson structure $\pi_{\bar{f}}$.  
\end{proof}


In both the constructions in the proofs of Theorem \ref{thm:TCvol-pres-flow-S4} and Corollary \ref{cor:TCvol-pres-flow-4-mfds}, the flows vanish in large portions of the ambient manifold. 
Only manifolds with zero Euler characteristic admit nonsingular flows.
Consequently, a general result about nonsingular flows that encompasses all smooth, closed, oriented $4$-manifolds is unattainable, as singularities would be inevitable for manifolds with nonzero Euler characteristic. 
It may nevertheless be interesting to wonder if such a non-singular, volume preserving and Turing complete flow can be constructed on smooth $4$-manifolds with vanishing Euler characteristic. 

 
 \medskip 
 
 The flow of a solution to a system like the one here only simulates a finite number of Turing machine steps (see \cite{CMPP21}). This phenomenon is related to the {\it Space-bounded Church Turing thesis}, introduced by Braverman, Schneider, and Rojas \cite{BSR15}. Storage space complexity and memory availability are studied as constraints for Turing machines modelled on compact spaces.  Moreover, {\it universal} Turing machines cannot be robustly simulated on a compact space \cite{BGH13}. 


 \medskip 
  
  Normally one wants to relate a volume preserving flow like those presented here to a solution of Euler's equations.  
   Our arguments rest on a background Poisson structure. 
   We have not assumed the existence of a Riemannian metric.  
  Hence, a model of hydrodynamics---perhaps some generalised Euler type equations---adapted to this particular situation would need to be specified first.
  


\subsubsection*{Acknowledgments} I am grateful to Misael Avenda\~no-Camacho for reading a previous version and suggesting improvements. Thanks to Jos\'e C. Ru\'iz Pantale\'on for useful conversations and correspondance about Poisson and Hamiltonian flows.



%
%

\providecommand{\bysame}{\leavevmode\hbox to3em{\hrulefill}\thinspace}

\end{document}